# Constrained Dynamic Control Allocation in the Presence of Singularity and Infeasible Solutions

David Buzorgnia and Ali Khaki-Sedigh

*Abstract*—Reliable controllers with high flexibility and performance are necessary for the control of intricate, advanced, and expensive systems such as aircrafts, marine vessels, automotive vehicles, and satellites. Meanwhile, control allocation has an important role in the control system design strategies of such complex plants. Although there are many proposed control allocation methodologies, few papers deal with the problems of infeasible solutions or system matrix singularity. In this paper, a pseudo inverse based method is employed and modified by the null space, least squares, and singular value decomposition concepts to handle such situations. The proposed method could successfully give an appropriate solution in both the feasible and infeasible sections in the presence of singularity. The analytical approach guarantees the solution with pre-defined computational burden which is a noticeable privilege than the linear and quadratic optimization methods. Furthermore, the algorithm complexity is proportionately grown with the feasible, infeasible, and singularity conditions. Simulation results are used to show the effectiveness of the proposed methodology.

## I. INTRODUCTION

Many advanced systems such as modern aircrafts, marine vessels, and automobiles have redundant actuators that can be used to enhance the closed loop performance, flexibility, and robustness. Control allocation is not a new subject in the control literature [1, 2]. However, the main developments in the control allocation techniques are related to the two last decades by demanding intricate systems and using computer instead of mechanical systems in aviation industry. Control allocation can be seen as an optimization problem with nonlinear limitations.

The earlier methods were based on linear algebra. Pseudo inverse is an appealing simple approach with minimum complexity [1, 3]. However, in the face of limitations, it cannot provide a practical solution. Daisy chaining methods and redistributed pseudo inverses (cascading generalized inverses) are the methods that try to solve the limitations problem, with partial success [1, 4]. Null Space intersection is another method that uses the null space to modify the pseudo inverse method. The large computational burden is its main drawback [1]. Geometric based methods are another alternative approach [1, 4]. Furthermore, attempts have been made to modify the idea (Direct Allocation) in order to prepare it for practical applications [5, 6].

Practical deficiencies of the earlier algebraic and geometric approaches [3] were a motivation to drive researchers towards the use of linear and quadratic programming methods. Fixed point method is a simple algorithm which can solve the problem [7]. However, its convergence is slow and strongly case dependent [3]. Active set method is another algorithm which solves the problem and decreases the computational time [8]. Also, it can find the exact answer. Interior point algorithm is a method which provides an approximate solution by implementing the $l_1$ norm [9, 10]. In [11], the $l_\infty$ norm is used in order to reduce the algorithm sensitivity when an actuator failure occurs. Nonlinear, robust, and adaptive methodologies are the other fields of control allocation techniques [12-14]. Singularity is a problem which can cause loss of controllability by converting the convex problem into a non-convex nonlinear problem. In [15], a cost function is proposed to avoid the singularity.

Most researches in control allocation have focused on the feasible desired virtual control vector based on the input vector constrains. In [14, 16], a method is proposed based on the pseudo inverse approach with reasonable computational burden. This method is revisited in this paper to further develop an algorithm applicable to both the feasible and infeasible cases. Furthermore, a practical method is proposed to solve the control allocation problem in the presence of singularity. Pseudo inverse, null space and least squares are the main tools that constructed the proposed method. The paper is organized as follows. In section 2, the linear control allocation problem is introduced as a convex optimization problem. The proposed method is presented in section 3. Section 4 provides the algorithm of the proposed method. In sections 5 and 6, a numerical example and simulation results demonstrate the applicability of the proposed method.

## II. PROBLEM STATEMENT

Consider a linear system with redundant inputs and the following state space equations:

$$\dot{x} = Ax + Bu, \qquad (1)$$
$$y = Cx, \qquad (2)$$

David Buzorgnia was with APAC Research Group, K. N. Toosi University of technology, Tehran, Iran (e-mail: d.b.buzorgnia@ieee.org)

Ali Khaki-Sedigh is with Center of Excellence in Industrial Control, Department of Electrical Engineering, K. N. Toosi University of Technology, Tehran, Iran (e-mail: sedigh@kntu.ac.ir)

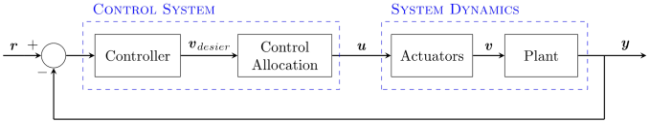

Fig. 1. Structure of a conventional linear system with controller and control allocation

where $A \in \mathbb{R}^{n \times n}$ is the system matrix, $B \in \mathbb{R}^{n \times m}$ is the input matrix, $C \in \mathbb{R}^{p \times n}$ is the output matrix, $x \in \mathbb{R}^n$ is the state vector, $u \in \Re \subset \mathbb{R}^m$ is the input vector and $\Re$ is the constrained control set due to physical limitations such as saturation, and $y \in \mathbb{R}^p$ is the output vector. Assume that $m > n$. As the number of inputs is more than the number of states, it is possible to control all the state variables with infinitely many solutions for the input vector corresponding to the same path of the state variables. Therefore, one can simultaneously reach the desired state variables and select one input vector from infinitely many candidates which satisfies other desired requirements. Without loss of generality, to simplify the problem, two separate design steps are considered. Let a new vector be defined as

$$v = Bu, \quad (3)$$

where $v \in \mathbb{R}^n$ is the virtual control vector. So, the system (1) can be written as follows

$$\dot{x} = Ax + v \quad (4)$$

Hence, the controller can find a unique solution for $v$ in (4). This solution is called the *desired virtual control vector* ($v_{desire}$). Assume that $v_{desire}$ is known. The problem is now to find a $u$ that satisfies

$$v_{desire} = Bu \quad (5)$$

With this new definition, the controller part is separated from the allocation part (fig. 1). After finding $v_{desire}$, (5) can be used to find $u$. That is, the control allocation process must find a $u$ to satisfy $v = v_{desire}$. Since there are infinitely many solutions for (5), a cost function is selected to customize the input vector. So, the problem is converted into an optimization problem. However, due to physical limitations, there are input vector constrains such as rate, maximum, and minimum saturations. Therefore, the control allocation problem is introduced as a convex optimization problem. In the case of non-convex problems, it is possible to break the problem into several convex problems. However, the global minimum cannot be guaranteed [17].

***Definition 1. (Control Allocation)*** Consider the linear system described by (1) and equivalently by (3) and (4). Assume that $m > n$ and $v_{desire}$ is known by the controller. Then, the input vector is determined by the following optimization problem:

$$\begin{aligned} &\text{minimize} \ f(u) \\ &\text{subject to} \ f_{limitations}(u) \leq 0, \quad (6)\\ &\quad Bu = v_{desire} \end{aligned}$$

where $f(u)$ is the cost function and $f_{limitations}(u)$ represents the physical limitations. Note that the above optimization problem must satisfy all the conditions for the convex optimization problem as given in [18].

***Definition 2. (Feasible and Infeasible Solutions)*** If there is at least one $u$ that simultaneously satisfies (5) and the physical limitations, the solution is feasible otherwise the solution is infeasible.

*Notations:* In this paper, scalars, vectors, matrices, and set of numbers are represented by $x$, $\boldsymbol{x}$, $X$, and $\mathbb{X}$ respectively. A scalar in subscript of a vector or matrix such as $x_i$ or $X_i$ represents the $i$ th element of the vector or $i$ th column (row) of the matrix if the number of columns is greater than the rows (the number of rows is greater than the columns). Furthermore, a set in subscript of a vector or matrix such as $x_S$ or $X_S$ demonstrates the corresponding elements of the vector or columns (rows) of the matrix if the number of columns is greater than the rows (the number of rows is greater than the columns). Through this paper, S represents the set of saturated elements with $k$ members and F is the set of free elements. Moreover, $i$ and $j$ are the representative of a saturated and free element, respectively.

## III. CONTROL ALLOCATION: THE PROPOSED METHOD

This section represents a new control allocation methodology based on the pseudo inverse, null space, and least squares concepts. The approach is to modify the pseudo inverse method in order to obtain an exact solution to the control allocation problem within the control constrains. The basic idea was previously introduced in [14, 16]. A quick review of the fundamental theory is presented in subsection 3.1. Subsection 3.2 demonstrates the developed feasible solution by considering a dynamic approach. Moreover, it provides a proper background for connecting to the infeasible section. Subsection 3.3 gives a practical solution for the infeasible desired virtual control vector in the presence of singularity.

### A. A review of the basic principles

In the absence of limitations, there exist infinitely many $u$ satisfying (5). One such $u$ is given by the pseudo inverse matrix, which minimizes the $l_2$ norm. However, in the face of control limitations, the derived solutions may not be implementable. Assume that $u$ is an arbitrary vector which satisfies (5). The aim is to modify it by a second vector as

$$u_{modified} = u - u_{nullity} \quad (7)$$

where $u_{nullity} \in \mathbb{R}^m$ is the correcting vector and $u_{modified}$ is the modified input vector which still satisfies (5). Multiplying (7) by $B$ gives

$$Bu_{modified} = Bu - Bu_{nullity} \quad (8)$$

so

$$v_{desire} = v_{desire} - Bu_{nullity} \quad (9)$$

hence

$$Bu_{nullity} = 0 \quad (10)$$

Therefore, the correcting vector $u_{nullity}$ must belong to the null space of $B$. So, (7) can be rewritten as follows

$$\boldsymbol{u}_{modified} = B^{-\dagger}\boldsymbol{v}_{desire} - N\boldsymbol{v}_{free} = \begin{bmatrix} B^{-\dagger} & -N \end{bmatrix} \begin{bmatrix} \boldsymbol{v}_{desire} \\ \boldsymbol{v}_{free} \end{bmatrix} \quad (11)$$

where $\boldsymbol{v}_{free} \in \mathbb{R}^{m-n}$ is a free virtual vector that is specified by the designer, $N$ is the null space matrix of $B$ and $B^{-\dagger}$ is the matrix pseudo inverse of $B$.

**Remark 1.** Equation (11) can be considered as a linear one to one function from $\mathbb{R}^m \to \mathbb{R}^m$. Hence, $\boldsymbol{v}_{free}$ maps $m-n$ elements of $\boldsymbol{u}_{modified}$. In other words, $\boldsymbol{v}_{free}$ can determine maximum $m-n$ elements of $\boldsymbol{u}_{modified}$.

**Lemma 1.** Assume that $\boldsymbol{u}$ satisfies (5) and be the pseudo inverse solution. Consider S as the set of saturated elements with $k$ members ($k \leq m-n$). To modify $\boldsymbol{u}_S$ in order to become equal to a desired vector ($\boldsymbol{u}_{desire}$) such that the modified input vector ($\boldsymbol{u}_{modified}$) still satisfies (5) with the minimum $l_2$ norm, the following equation could be used:

$$\boldsymbol{u}_{modified} = \boldsymbol{u} - N\boldsymbol{v}_{free} \quad (12)$$

where $\boldsymbol{v}_{free}$ is

$$\boldsymbol{v}_{free} = N_S^T (N_S N_S^T)^{-1} \Delta \quad (13)$$

and $T$ is the transpose symbol. $\Delta$ is defined as

$$\Delta = \boldsymbol{u}_S - \boldsymbol{u}_{desire} \quad (14)$$

**Proof:** From (11), we have,

$$\|\boldsymbol{u}_{modified}\|_2 = \|\boldsymbol{u} - N\boldsymbol{v}_{free}\|_2 \quad (15)$$

so

$$\|\boldsymbol{u}_{modified}\|_2 \leq \|\boldsymbol{u}\|_2 + \|N\boldsymbol{v}_{free}\|_2 \quad (16)$$

where $\|\cdot\|_2$ denotes the $l_2$ norm. As $\boldsymbol{u}$ is known from the pseudo inverse, $\|N\boldsymbol{v}_{free}\|_2$ must be minimized. Assume that $N$ is known, so $\|\boldsymbol{v}_{free}\|_2$ must be minimized. The goal is to correct just the $k$ elements of $\boldsymbol{u}$ which are determined by the saturated set and the other elements are free to move in order to obtain the minimum $l_2$ norm of the modified input vector. So

$$N_S \boldsymbol{v}_{free} = \Delta \quad (17)$$

The pseudo inverse solution could be used to obtain $\boldsymbol{v}_{free}$ with the possible minimum $l_2$ norm. Minimizing $\|\boldsymbol{v}_{free}\|_2$ could guarantee the minimum $l_2$ norm of $\boldsymbol{u}_{modified}$. ∎

Lemma 1 is based on the assumption that $N$ is known. From (13), it is obvious that $\|\boldsymbol{v}_{free}\|_2$ depends on $N$.

**Lemma 2.** To minimize $\|\boldsymbol{v}_{free}\|_2$ in (13), select $N$ such that

$$\underline{\sigma} = \overline{\sigma} \quad (18)$$

where $\underline{\sigma}$ and $\overline{\sigma}$ are the minimum and maximum singular values of $N$, respectively.

**Proof:** The $l_2$ norm of (17) gives

$$\|N_S \boldsymbol{v}_{free}\|_2 = \|\Delta\|_2 \quad (19)$$

so

$$\|\boldsymbol{v}_{free}\|_2 = \frac{1}{\sigma'} \|\Delta\|_2 \quad (20)$$

where $\sigma'$ is a gain between the maximum ($\overline{\sigma}'$) and minimum ($\underline{\sigma}'$) singular values of $N_S$ that gives the equality. The object of Lemma 1 is to minimize $\|N\boldsymbol{v}_{free}\|_2$. So,

$$\|N\boldsymbol{v}_{free}\|_2 = \sigma \|\boldsymbol{v}_{free}\|_2 \quad (21)$$

where $\sigma$ is a gain between the maximum and minimum singular values of $N$ that results in equality. Equation (20) and (21) gives

$$\|N\boldsymbol{v}_{free}\|_2 = \frac{\sigma}{\sigma'} \|\Delta\|_2 \quad (22)$$

The term $\Delta$ is known and therefore $\sigma/\sigma'$ must be minimized. From Appendix A, Lemma 11, we know that $\sigma' \leq \sigma$. Consider the worst case that $\sigma = \overline{\sigma}$ and $\sigma' = \underline{\sigma}'$. In order to minimize $\sigma/\sigma'$, $\overline{\sigma}$ should be decreased and $\underline{\sigma}'$ should be increased as much as possible. We do not have a direct access to $\underline{\sigma}'$ since $N_S$ depends on $N$. From lemma 11, the minimum amount for $\overline{\sigma}$ is $\underline{\sigma}$. ∎

The singular value decomposition of the matrix $B$ could be written as $B = U \Sigma V^T$ where $U$ and $V$ are the orthonormalized eigenvectors of $BB^T$ and $B^T B$ respectively and the diagonal elements of $\Sigma$ are the singular values. The last $m-n$ columns of $V$ constitute the null space matrix of $B$ with the singular values equal to one.

### B. Feasible Desired Virtual Control Vector

In this paper, a criterion is proposed to provide a comparison assessment of the elements of the input vector. This comparative criterion is based on the relative positions of the corresponding maximum and minimum limitations.

**Definition 3. (Weighted Distance Vector)** Consider the input vector $\boldsymbol{u}$ with the physical upper and lower limitations of $\boldsymbol{u}_{min}$ and $\boldsymbol{u}_{max}$, respectively. Define the bandwidth center vector $\boldsymbol{u}_{center}$ as

$$\boldsymbol{u}_{center} = \frac{\boldsymbol{u}_{max} - \boldsymbol{u}_{min}}{2} \quad (23)$$

and the border vector as $\boldsymbol{u}_{border} = [\boldsymbol{u}_{border_l}]_{l=1}^{m}$ where

$$\boldsymbol{u}_{border_l} = \begin{cases} \boldsymbol{u}_{max_l} & \text{if } \boldsymbol{u}_l \geq \boldsymbol{u}_{center_l} \\ \boldsymbol{u}_{min_l} & \text{if } \boldsymbol{u}_l < \boldsymbol{u}_{center_l} \end{cases} \quad (24)$$

Now, the weighted distance vector can be defined as $\boldsymbol{w} = [\boldsymbol{w}_l]_{l=1}^{m}$ where

$$\boldsymbol{w}_l = \frac{\boldsymbol{u}_l - \boldsymbol{u}_{center_l}}{\boldsymbol{u}_{border_l} - \boldsymbol{u}_{center_l}} \quad (25)$$

**Definition 4. (Saturated Set)** By considering definition 3,

the saturated set is described as
$$S = \{l : l = 1, ..., m \mid w_l = \|w\|_\infty\} \quad (26)$$
where $l$ is the element number corresponding to $u_l$. In other words, the indices of the elements of the input vector which belong to the infinity norm of the weighted distance vector constitute the members of the saturated set.

There are a few solving strategies in the control allocation methods. Saving the direction of the desired virtual control vector, for instance, is the main criterion of the direct allocation method proposed in [1]. A tradeoff between preserving the minimum $\|u\|_2$ and $\|v_{desire} - v\|_2$ is the criterion of the linear and quadratic programing in [3]. The programming method proposed in [11] minimizes $\|u\|_1$ or $\|u\|_\infty$ instead of $\|u\|_2$. The criterion in this paper is to maintain the minimum $\|u\|_2$ while $\|w\|_\infty$ is being reduced.

**Lemma 3.** Consider
$$\delta = |\Delta_t| \quad (27)$$
where $\Delta$ comes from (14) and $t$ is a specific member of the saturated set. Computing $\Delta$ with
$$\Delta = \bar{\Delta}\delta \quad (28)$$
where $\bar{\Delta} = [\bar{\Delta}_i]_{\forall i \in S}$ and
$$\bar{\Delta}_i = \frac{u_{border_i} - u_{center_i}}{|u_{border_t} - u_{center_t}|} \quad (29)$$
Applying the above $\Delta$ in (13) reduces $\|w_{modified}\|_\infty$ in (12).

*Proof:* As $u_S$ detonates the elements corresponding to the saturated set, it is obvious that $w_S = \|w\|_\infty$. To reduce $\|w\|_\infty$, $w_i$ for all the saturated elements should be decreased with the same proportion, which requires $w_{modified_i} = w_{modified_t}$ $\forall i \in S$. Assume that $u_{border_i}$ is the same before and after modifications. So, by considering $w_S$ before and after modifications
$$w_i - w_{modified_i} = w_t - w_{modified_t} \quad (30)$$
and using (14), it can be simplified as
$$\frac{\Delta_i}{u_{border_i} - u_{center_i}} = \frac{\Delta_t}{u_{border_t} - u_{center_t}} \quad (31)$$
It could be easily shown that the sign of the numerator and denominator of the fractions of (31) are the same. So, taking the absolute value of (31) yields
$$\frac{\Delta_i}{u_{border_i} - u_{center_i}} = \frac{|\Delta_t|}{|u_{border_t} - u_{center_t}|} \quad (32)$$
Hence, (32) gives the opportunity to compute $\Delta$ based on a specific element ($t$ th element). ∎

**Remark 2.** By using lemma 3, lemma 1 could be written based on the scalar $\delta$ as follows
$$u_{modified} = u - u_{reduction}\delta \quad (33)$$
where $u_{reduction}$ is
$$u_{reduction} = NN_S^T(N_S N_S^T)^{-1}\bar{\Delta} \quad (34)$$

**Definition 5. (Intersection Value)** The intersection value of the $j$ th element of the modified input vector (33) that $j \in F$, is a specific value of $\delta$ which results in
$$w_{modified_j} = w_{modified_t} \quad (35)$$
Note that the number of the intersection values is equal to the number of the free elements ($m - k$). To compute the intersection value corresponding to $j$ th element, (25), (33), and (35) are used as follows
$$\frac{u_j - u_{reduction_j}\delta_j - u_{center_j}}{u_{border_j} - u_{center_j}} = \frac{u_t - u_{reduction_t}\delta_j - u_{center_t}}{u_{border_t} - u_{center_t}} \quad (36)$$
where $\delta_j$ is the intersection value of $j$ th element. By a simple manipulation, (36) yields
$$\delta_j = \left(\frac{u_{reduction_t}}{u_{border_t} - u_{center_t}} - \frac{u_{reduction_j}}{u_{border_j} - u_{center_j}}\right)^{-1} \times \left(\frac{u_t - u_{center_t}}{u_{border_t} - u_{center_t}} - \frac{u_j - u_{center_j}}{u_{border_j} - u_{center_j}}\right) \quad (37)$$

**Remark 3.** By Considering $u_{border_j} = u_{max_j}$
$$u_{border_j} - u_{center_j} = \frac{u_{max_j} + u_{min_j}}{2} \quad (38)$$
Consecutively, if $u_{border_j} = u_{min_j}$, then
$$u_{border_j} - u_{center_j} = -\frac{u_{max_j} + u_{min_j}}{2} \quad (39)$$

**Lemma 4.** The intersection value $\delta_j$ which is computed by (37) must satisfy the following conditions:
- $\delta_j > 0$
- $\dfrac{u_t - u_{reduction_t}\delta_j - u_{center_t}}{u_{border_t} - u_{center_t}} > 0$

If $\delta_j$ violates any of the above conditions, let $u_{border_j} - u_{center_j}$ as $-(u_{border_j} - u_{center_j})$ and re-compute (37).

*Proof:* From lemma 3, $\delta_j$ is defined as a positive scalar (the first condition). Moreover, it is assumed that $u_{border_S}$ is the same for the non-modified and the modified input vector (the second condition). Let $\delta_j$ be negative (violating the first condition). From (37), the term
$$\left(\frac{u_t - u_{center_t}}{u_{border_t} - u_{center_t}} - \frac{u_j - u_{center_j}}{u_{border_j} - u_{center_j}}\right) \geq 0 \quad (40)$$
is always positive since it corresponds to $w_t - w_j$. Therefore, the term
$$\left(\frac{u_{reduction_t}}{u_{border_t} - u_{center_t}} - \frac{u_{reduction_j}}{u_{border_j} - u_{center_j}}\right) < 0 \quad (41)$$
must be negative. So, from (41) it yields
$$\frac{u_{border_j} - u_{center_j}}{u_{border_t} - u_{center_t}} < \frac{u_{reduction_j}}{u_{reduction_t}} \quad (42)$$

Without loss of generality, equation (42) is obtained from (41) based on the assumption that the direction of inequality remains the same. $u_{border_t}$, $u_{reduction_t}$, and $u_{reduction_j}$ are fixed during the modifications. However, $u_j$ could pass over $u_{center_j}$ which results in changing $u_{border_j}$. Based on remark 3, changing $u_{border_j}$ gives

$$\frac{-(u_{border_j} - u_{center_j})}{u_{border_t} - u_{center_t}} < \frac{u_{reduction_j}}{u_{reduction_t}} \quad (43)$$

A negative sign alters the inequality direction and consequently makes $\delta_j > 0$. Similarly, let condition two be violated. Therefore, the left and right signs of (36) are not equal. By changing $u_{border_j}$ and rewriting (36), it yields

$$\frac{u_j - u_{reduction_j}\delta_j - u_{center_j}}{-(u_{border_j} - u_{center_j})} = \frac{u_t - u_{reduction_t}\delta_j - u_{center_t}}{u_{border_t} - u_{center_t}} \quad (44)$$

where the produced negative sign makes the equality true. ∎

**Lemma 5.** The specific $\delta$ which gives the following relation

$$w_{modified_j} = w_{modified_i} = \|w_{modified}\|_\infty \quad (45)$$

is obtained by

$$\delta = \min\left(\delta_j \mid \forall j, j \in F\right) \quad (46)$$

*Proof:* By increasing $\delta$, lemma 1 and 3 show that $w_{modified_S}$ would be decreased. Moreover, increasing $\delta$ leads to an increase in $\|u_{modified}\|_2$. So, it can be concluded that increasing $\delta$ also increases $\|w_{modified_F}\|_\infty$. Assume that $j$ is the element that holds (45). Further an increase in $\delta_j$ leads to $w_{modified_j} > w_{modified_i}$ and consequently $\|w_{modified}\|_\infty > w_{modified_i}$. The minimum $\delta_j$ guarantees that no free element satisfies $w_{modified_j} > w_{modified_i}$. ∎

**Theorem 1. (Control Allocation for Feasible Desired Virtual Control Vector)** Let $v_{desire}$ be given and feasible. To obtain $u$ in order to have the minimum $l_2$ norm and satisfy the limitations, compute $u$ by the pseudo inverse method. If it lies within the limitations, it would be the answer. If not, compute the intersection values for all the free elements by (37) with considering lemma 4. Select the appropriate $\delta$ among the intersection values with the criterion introduced by lemma 5. If the selected $\delta$ gives

$$\frac{u_t - u_{reduction_t}\delta - u_{center_t}}{u_{border_t} - u_{center_t}} > 1 \quad (47)$$

then modify the input vector by remark 2 and repeat the procedure from updating the saturated set with definition 4 and finding the intersection value. Else, re-compute $\delta$ with

$$\delta = |u_t - u_{border_t}| \quad (48)$$

and modify the input vector by remark 2.

**Remark 4.** Pseudo inverse method gives the minimum $l_2$ norm of the input vector. Modifying input vector by remark 2, which is based on lemma 1, 2, and 3, preserves the minimum $l_2$ norm while reducing $\|w\|_\infty$. Lemma 4 and 5 provide a reliable way for obtaining the proper $\delta$. Since the whole procedure is based on modifying the input vector by (7), repeating the algorithm and correcting the input vector with several $u_{nullity}$ could still preserve the minimum $l_2$ norm of the input vector. Note that reducing $\|w\|_\infty$ increases $l_2$ norm of the input vector. To obtain the minimum $l_2$ norm, (48) is used to prevent $\|w\|_\infty$ from unnecessary reduction which causes $\|w\|_\infty = 1$.

### C. Infeasible Desired Virtual Control Vector

The following lemmas and theorems are based on the assumption that $v_{desire}$ is infeasible.

**Lemma 6.** To modify the input vector in order to reduce $\|w\|_\infty$ while preserving the minimum $\|v_{desire} - v\|_2$, the following relation could be used.

$$\begin{bmatrix} u_{modified_S} \\ u_{modified_F} \end{bmatrix} = \begin{bmatrix} u_S \\ 0 \end{bmatrix} + \begin{bmatrix} -\overline{\Delta} & 0 \\ B_F^\dagger B_S \overline{\Delta} & B_F^\dagger \end{bmatrix} \begin{bmatrix} \delta \\ v_{desire} - B_S u_S \end{bmatrix} \quad (49)$$

where $B_F^\dagger$ is the matrix inverse of $B_F$.

*Proof:* Lemma 3 gives a $\Delta$ which could decrease $\|w\|_\infty$ with a simple computation. The modified saturated elements could be rewritten as

$$u_{modified_S} = u_S - \overline{\Delta}\delta \quad (50)$$

By separating the free and saturated elements, it yields

$$v_{desire} = B_S u_{modified_S} + B_F u_{modified_F} \quad (51)$$

Assume that $B_F^\dagger$ is the inverse of $B_F$ which minimizes $\|v_{desire} - v\|_2$. So,

$$u_{modified_F} = B_F^\dagger \left(v_{desire} - B_S u_{modified_S}\right) \quad (52)$$

∎

Consider $\overline{B} \in \mathbb{R}^{n \times p}$ as a singular sub-matrix of $B$ with $p \leq m$. There are many iterative methods to estimate the matrix inverse of a singular matrix which are dependent on the matrix structure [19, 20]. For control allocation, there are two important factors in singular matrix inversion:

- The product must be close to the identity matrix.
- The elements of the matrix inverse should prevent saturation of the input vector.

The following theorem offers a solution by considering the stated criteria.

**Theorem 2.** Assume that $\overline{B} \in \mathbb{R}^{n \times p}$ is a singular matrix with $rank(\overline{B}) = r$ where $r < n$. Let $\overline{B} = UDV^*$ be the corresponding singular value decomposition, then define $B_{augmented}$ as follows

$$B_{augmented} = \begin{bmatrix} \overline{B} & U_{r:n} \end{bmatrix} \quad (53)$$

where $r:n$ denotes the last $n-r$ columns of $U$. The pseudo inverse of $B_{augmented}$ is $P = B_{augmented}^T \left( B_{augmented} B_{augmented}^T \right)^{-1}$ and $\bar{B}^\dagger = P_{1:m}$ is the matrix inverse of $\bar{B}$ where $1:m$ denotes the first $m$ rows of $P$.

*Proof:* $U$ is the left-singular vectors of $\bar{B}$. Since the last $r-n$ eigenvalues of $\bar{B}$ are zero, the corresponding columns in $U$ consist of the n-dimensional null space vectors. By adding this vectors to $\bar{B}$, $B_{augmented}$ would become full rank ($rank(B_{augmented}) = n$). Pseudo inverse guarantees the matrix inverse with the minimum $l_2$ norm. The last $r$ rows of $P$ which correspond to the null space vectors should be eliminated. By obtaining the minimum $l_2$ norm of the input vector, input saturation occurrence can be largely avoided. ∎

*Proposition 1.* For a matrix $B$ with rank $n$ if there exist some $\bar{B} \subset B$ with $rank(\bar{B}) < n$, that is there are linearly dependent columns in $B$, pseudo invers could still successfully compute the matrix inverse of $B$.

*Proposition 2. (Infeasible Conditions).* As stated in definition 2, if there is no $u$ that simultaneously satisfies (5) and the physical constrains, the solution is infeasible. Infeasibility occurs in the following cases:

**Condition 1.** If the number of the saturated elements is larger than the number of the free virtual vector elements ($k > (m-n)$).

**Condition 2.** If the rank of $B$ is lower than $n$.

**Condition 3.** If $\det(N_S N_S^T) = 0$

*Proof:* The first condition can be proved by remark 1. The second condition is obvious and dependent rows or columns in $B$ make $\det(BB^T) = 0$. Hence, $B$ does not have a matrix inverse and the solution is totally infeasible. The last condition occurs when there exist some dependent columns in $B$ and consequently some dependent rows in $N$. In such cases, the solution is infeasible since (13) has no answer. Note that the two last conditions are the results of a singular matrix. ∎

*Lemma 7.* Based on the dimension and rank of $B_F$, computing $B_F^\dagger$ by one of the following cases could give the minimum $\|v_{desire} - v\|_2$ in lemma 6.

**Case 1.** If $rank(B_F) = n$, compute $B_F^\dagger$ with the pseudo inverse method.

**Case 2.** If $rank(B_F) = m-k$, compute $B_F^\dagger$ with the least squares method.

**Case 3.** If $rank(B_F) < m-k$, compute $B_F^\dagger$ with theorem 2.

*Proof:* In case 1, $B_F$ is a full rank matrix. Based on proposition 1, pseudo inverse is the best candidate to obtain the matrix inverse with minimum $\|u_{modified_F}\|_2$. In case 2, rank of $B_F$ is equal to the number of the free elements and is lower than $n$. Least squares can obtain $u_{modified_F}$ to minimize $\|v_{desire} - v\|_2$. In case 3, $B_F$ is rank deficient. Theorem 2 provides a fast and reliable solution to compute the pseudo inverse for the singular matrices. ∎

Since the modification strategies of the feasible and infeasible sections are reducing the infinity norm of the Weighted Distance Vector introduced by definitions 3 and 4, the similar description of Intersection Value stated by definition 5 could be implemented here. To compute the intersection value corresponding to the $j$ th element, (25), (35), and (49) are implemented as follows

$$\frac{B_j^\dagger \left( v_{desire} - B_S u_S \right) + B_j^\dagger B_S \bar{\Delta} \delta_j - u_{center_j}}{u_{border_j} - u_{center_j}} = \frac{u_t - \bar{\Delta}_t \delta_j - u_{center_t}}{u_{border_t} - u_{center_t}} \quad (54)$$

where $\delta_j$ is the intersection value corresponding to the $j$ th element. By a simple manipulation, it yields

$$\delta_j = \left( \frac{\bar{\Delta}_t}{u_{border_t} - u_{center_t}} + \frac{B_j^\dagger B_S \bar{\Delta}}{u_{border_j} - u_{center_j}} \right)^{-1} \times \left( \frac{u_t - u_{center_t}}{u_{border_t} - u_{center_t}} - \frac{B_j^\dagger (v_{desire} - B_S u_S) - u_{center_j}}{u_{border_j} - u_{center_j}} \right) \quad (55)$$

*Lemma 8.* The intersection value $\delta_j$ satisfies the following inequalities:

- $\delta_j > 0$

- $\dfrac{u_t - \bar{\Delta}_t \delta_j - u_{center_t}}{u_{border_t} - u_{center_t}} > 0$

If $\delta_j$ does not satisfy any of the above inequalities, consider $u_{border_j} - u_{center_j}$ as $-(u_{border_j} - u_{center_j})$ and re-compute (55).

*Proof:* The explanation for the necessity of the two conditions is the same as stated in lemma 4. Assume that

$$\frac{B_j^\dagger (v_{desire} - B_S u_S) - u_{center_j}}{u_{border_j} - u_{center_j}} \leq \frac{u_t - u_{center_t}}{u_{border_t} - u_{center_t}} \quad (56)$$

So, the only term in (55) that makes $\delta_j$ negative is

$$\left( \frac{\bar{\Delta}_t}{u_{border_t} - u_{center_t}} + \frac{B_j^\dagger B_S \bar{\Delta}}{u_{border_j} - u_{center_j}} \right) < 0 \quad (57)$$

So

$$\frac{u_{border_j} - u_{center_j}}{u_{border_t} - u_{center_t}} < -\frac{B_j^\dagger B_S \bar{\Delta}}{\bar{\Delta}_t} \quad (58)$$

Without loss of generality, assume that the direction of inequality in (57) is the same as (58). Based on remark 3, changing the border of $u_{border_j}$ leads to

$$\frac{-(u_{border_j} - u_{center_j})}{u_{border_t} - u_{center_t}} < -\frac{B_j^\dagger B_S \bar{\Delta}}{\bar{\Delta}_t} \quad (59)$$

The negative sign inverses the inequality direction and consequently makes $\delta_j$ positive. Similarly, assume that the second condition is violated. Hence, $u_{border_t}$ is changed which causes a negative sign in (54). By changing $u_{border_j}$, (54) could be preserved. ∎

***Lemma 9.*** The criterion for selecting the appropriate $\delta$ is the same as stated in lemma 5.

***Proof:*** Increasing $\delta$ in (50), decreases $\|u_{modified_S}\|_2$ and $w_S$. Moreover, increasing $\delta$ in (52) gradually increases $\|w_F\|_\infty$. It is obvious that the minimum intersection value belongs to the first free element that gives $w_j = w_S$. Any further increase in $\delta$ leads to $w_j > w_S$ and consequently $w_j > \|w\|_\infty$. ∎

Lemma 8 and 9 are based on the assumption that (56) is true. However, in the case of singularity, (56) may not hold, since $u_{modified_F}$ is not computed based on $u_F$ (see (49)). Hence, it is possible for $\delta = 0$ in (49), $w_{modified_j} > w_{modified_t}$. In this case, the following strategy given in Lemma 10 can be used.

***Lemma 10.*** For $\delta = 0$ in (49) if any free element has $w_{modified_j} > w_{modified_t}$ (called the escaped element), select the element corresponding to the maximum $\delta'_j$ as the saturated element, where

$$\delta'_j = \frac{w_j - w_t}{u'_{reduction_j}}, \forall j \in F, w_j > w_t \tag{60}$$

$$u'_{reduction_j} = \frac{B_j B_S \overline{\Delta}}{u_{border_j} - u_{center_j}} \tag{61}$$

and modify it by

$$u_{modified_j} = w_t(u_{border_j} - u_{center_j}) + u_{center_j} \tag{62}$$

which gives $w_{modified_j} = w_t$.

***Proof:*** Consider $u_S$ is fixed ($u_{modified_S} = u_S$). By adding a negative sign to $\delta$ in (49), a reverse procedure is implemented in order to reduce $w_F$ until $w_{modified_j} \leq w_{modified_t}$ for all the escaped free elements. Equation (49) can be rewritten with the negative sign as

$$u_{modified_F} = B_F^\dagger(v_{desire} - B_S u_S) - B_F^\dagger B_S \overline{\Delta} \delta' \tag{63}$$

Assume that $j$ is one of the escaped free elements. The objective is to have $w_{modified_j} = w_t$. So

$$w_j - w_{modified_j} = w_j - w_t \tag{64}$$

By implementing (25) and (63), (64) could be written as

$$w_j - w_t = \frac{B_j^\dagger B_S \overline{\Delta} \delta'_j}{u_{border_j} - u_{center_j}} \tag{65}$$

Since $w_j$ and $w_t$ are known, $\delta'_j$ can be computed. Note that $B_j^\dagger$ is different for every free element which is caused to obtain a different $\delta'_j$ for every escaped free element. Selecting the element corresponding to the maximum $\delta'_j$ guarantees $w_{modified_j} \leq w_{modified_t}$ for all the escaped free elements. ∎

***Theorem 3. (Control Allocation for the Infeasible Desired Virtual Control Vector)*** Assume that $v_{desire}$ is infeasible (based on the conditions stated in proposition 2). To obtain a $u$ that satisfies the limitations and leads to a minimum $\|v_{desire} - v\|_2$ and reduces $\|w\|_\infty$, compute $B_F^\dagger$ by lemma 7. To ensure that the saturated elements are selected properly, use lemma 10. Obtain the intersection points for free elements by (55) and lemma 8 and choose the appropriate intersection point by lemma 9. Modify the input control vector by lemma 6. If the selected $\delta$ gives

$$\frac{u_t - \overline{\Delta}\delta - u_{center_t}}{u_{border_t} - u_{center_t}} > 1 \tag{66}$$

then repeat the above procedure. Else, re-compute $\delta$ by (48) and modify the input vector by lemma 6.

***Remark 5.*** If there is no free elements and all the elements belong to the saturated set, then reduce the modified input vector until

$$u_{modified} = u_{border} \tag{67}$$

To further clarify the main points, a flowchart of the proposed algorithm is depicted in fig. 2. Note that the left and right columns are associated to the feasible and infeasible sections respectively.

***Proposition 3.*** The maximum number of iterations in the feasible section (the loop in the left column of fig. 2) is equal to $m - n$ and the maximum number of iterations for both feasible and infeasible sections (the loops in the left and right columns of fig. 2) is equal to $m$. In other words, the maximum number of the iterations is known and is equal to the number of the elements of the input vector.

## IV. A NUMERICAL EXAMPLE

Consider $B$, $u_{min}$, $u_{max}$, and $v_{desired}$ as follows

$$B = \begin{bmatrix} 1 & 1 & 1 & 1 & 1 \\ 1 & 1 & 1 & 0 & 0 \\ 1 & 0 & 0 & 0 & 0 \end{bmatrix}$$

$$u_{min} = [-1 \quad 0.2 \quad -1 \quad -0.4 \quad -0.2]^T$$

$$u_{max} = [1.2 \quad 1 \quad 0 \quad 0.6 \quad 0.1]^T$$

$$v_{desire} = [1.4 \quad 1 \quad -1]^T$$

To find the input vector, follow the algorithm depicted in fig. 2. Note that $B$ is a full rank matrix with dependent columns. Compute the input vector and $w$ by the pseudo inverse method and definition 3.

$$u = [-1 \quad 1 \quad 1 \quad 0.2 \quad 0.2]^T$$

$$w = [1 \quad 1 \quad 3 \quad 0.2 \quad 1.667]^T$$

$\|w\|_\infty$ is greater than 1. Obtaining the saturated set by definition 4 gives $S = \{3\}$. Compute $\delta_j$ for all the free elements by (37) and lemma 4.

$$[\delta_1 \quad \delta_2 \quad \delta_4 \quad \delta_5] = [1 \quad 0.444 \quad 1.4 \quad 0.667]$$

Selecting the appropriate $\delta$ by lemma 5 gives $\delta = 0.444$ which belongs to the second element. Now, the input vector can be modified by (33) which gives

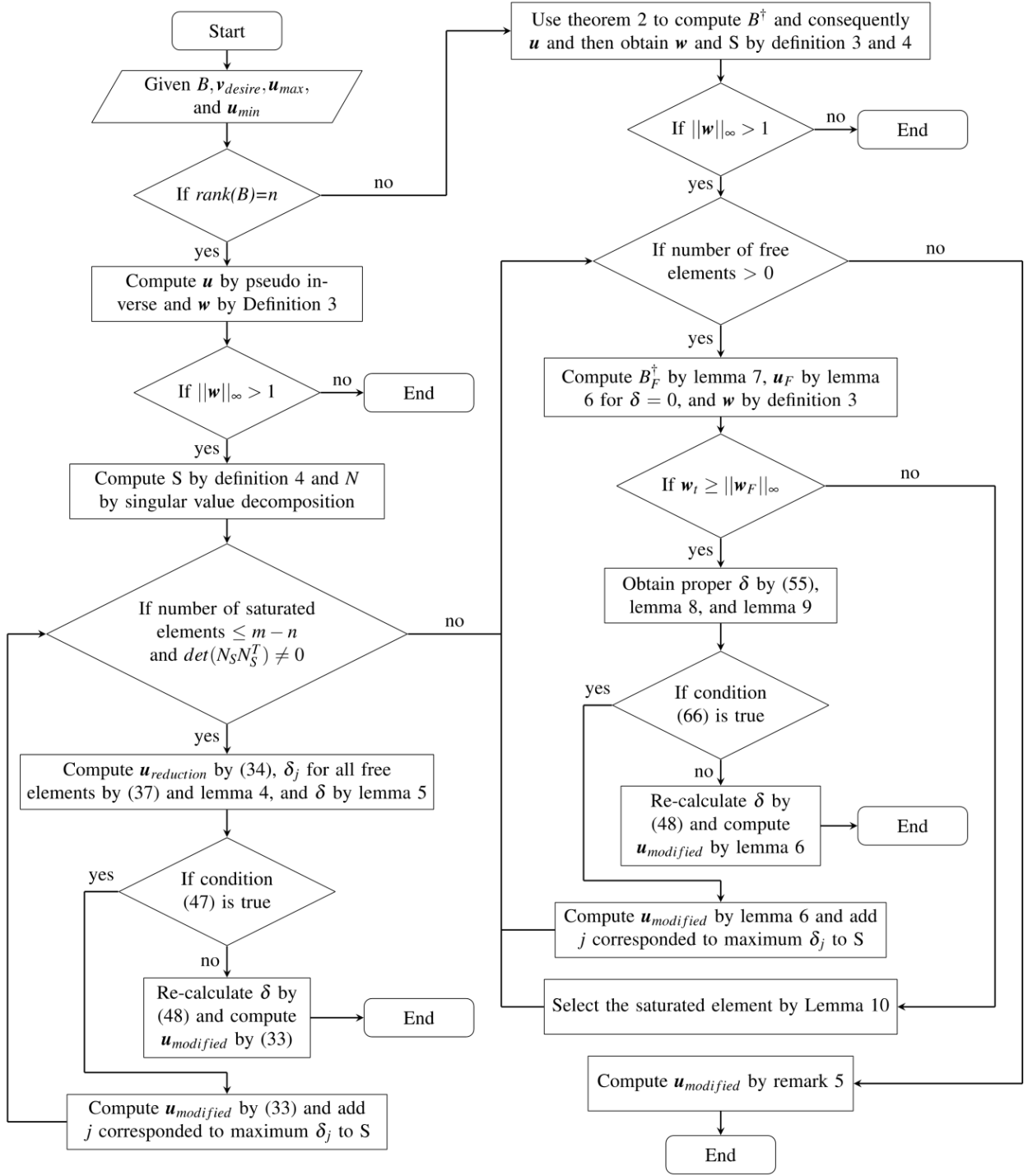

Fig. 2. Flowchart of the proposed method

$$\boldsymbol{u}_m^1 = \begin{bmatrix} -1 & 1.444 & 0.556 & 0.2 & 0.2 \end{bmatrix}^T$$

The superscript of $\boldsymbol{u}_m$ shows the number of iterations or modifications. The saturated set can now be written as $S = \{2 \ 3\}$. Since condition (47) is true, the above procedure should be repeated. Before continuing, the infeasible conditions stated in proposition 2 should be checked. Since columns 2 and 3 of $B$ are linearly dependent, $\det(N_S N_S^T)$ is zero and causes the infeasibility condition. So, compute $B_F^\dagger$ by lemma 7. Since $rank(B_F) < m - k$, $B_F^\dagger$ should be computed by theorem 2 which gives

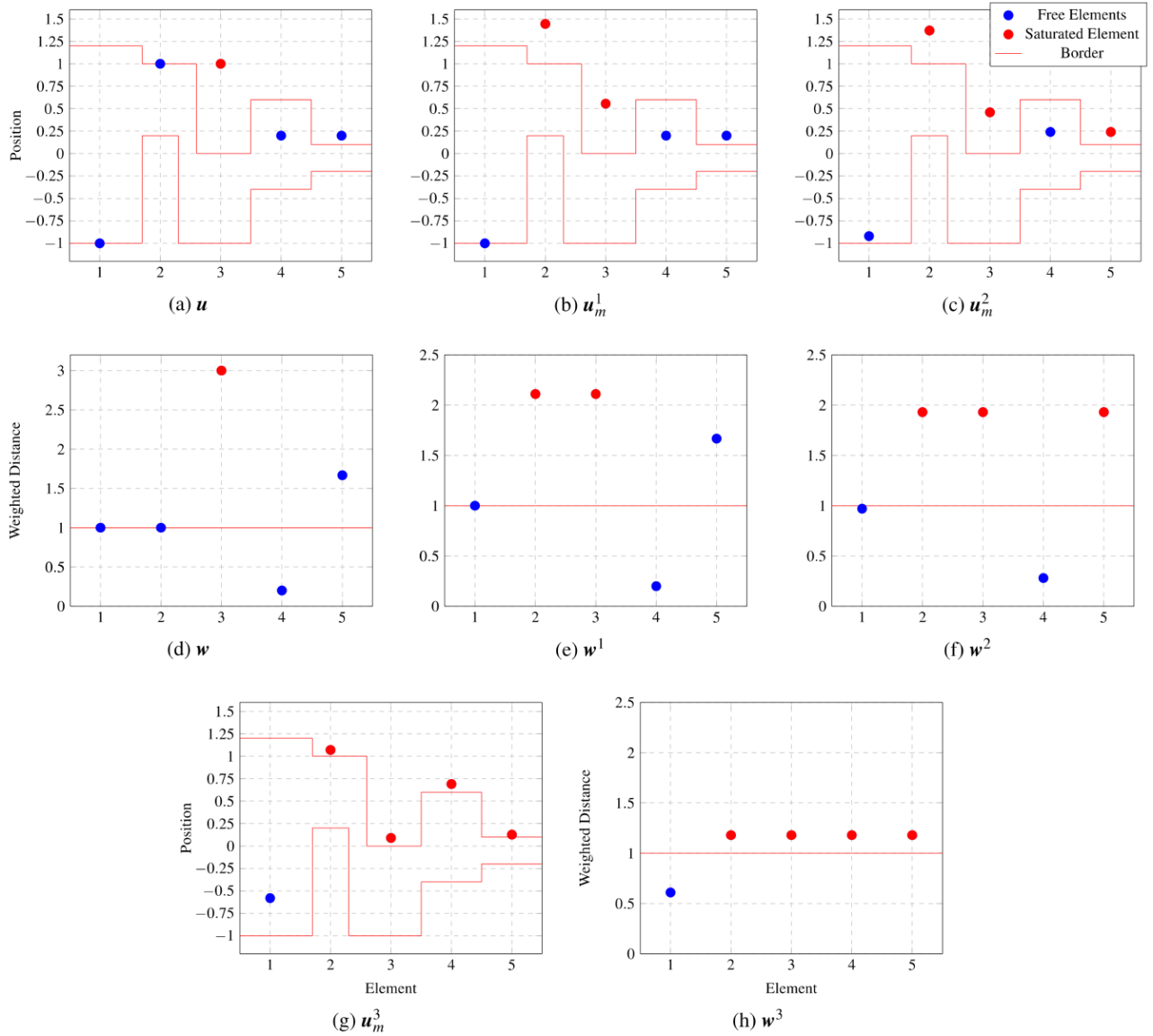

Fig. 3. Position of the elements of the input and weighted distance vectors during the modification

$$B_F^\dagger = \begin{bmatrix} 0 & 0.5 & 0.5 \\ 0.5 & -0.25 & -0.25 \\ 0.5 & -0.25 & -0.25 \end{bmatrix}$$

Computing $u_F$ for $\delta = 0$ by lemma 6 gives $u_F = [-1 \ 0.2 \ 0.2]^T$. So

$$w^1 = [1 \ 2.111 \ 2.111 \ 0.2 \ 1.667]$$

Calculating $\delta_j$ by (55) and lemma 8 gives

$$[\delta_1 \ \delta_4 \ \delta_5] = [0.940 \ 0.659 \ 0.088]$$

Selecting the proper $\delta_j$ by lemma 9 yields $\delta = 0.088$ which belongs to the fifth element. Since condition (66) is true, modify the input vector by lemma 6

$$u_{modified}^2 = [-0.92 \ 1.373 \ 0.467 \ 0.24 \ 0.24]^T$$

and update the saturated set $S = \{2 \ 3 \ 5\}$. Repeat the procedure and compute $B_F^\dagger$. Since $rank(B_F) = m - k$, $B_F^\dagger$ can be computed by the least squares which gives

$$B_F^\dagger = \begin{bmatrix} 0 & 0.5 & 0.5 \\ 1 & -0.5 & -0.5 \end{bmatrix}$$

Computing $u_F$ for $\delta = 0$ gives $u_F = [-0.92 \ 0.24]^T$. So

$$w^2 = [0.972 \ 1.933 \ 1.933 \ 0.28 \ 1.933]$$

Calculating $\delta_j$ gives

$$[\delta_1 \ \delta_4] = [0.851 \ 0.376]$$

Hence, $\delta = 0.376$ which belongs to the fourth element. Condition (66) is true. So

$$u_{modified}^3 = [-0.581 \ 1.072 \ 0.091 \ 0.691 \ 0.127]^T$$

and $S = \{2 \ 3 \ 4 \ 5\}$. Repeat the procedure and compute $B_F^\dagger$ by the least squares. $u_F$ for $\delta = 0$ is $u_F = [-0.581]$ and

$$\boldsymbol{w}^3 = [0.62 \quad 1.182 \quad 1.182 \quad 1.182 \quad 1.182]$$

and $\delta = 0.517$. Condition (66) is not true anymore. So, re-calculate $\delta$ with (48) which gives $\delta = 0.0909$ and the modified input vector is

$$\boldsymbol{u}^4_{modified} = [-0.433 \quad 1 \quad 0 \quad 0.6 \quad 0.1]^T$$

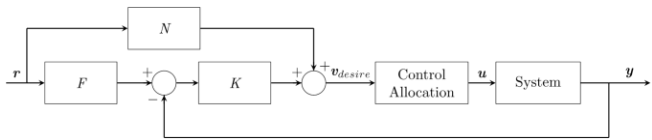

Fig. 4. Over-actuated Satellite Launch Vehicle (SLV) with linear quadratic set point control and control allocation

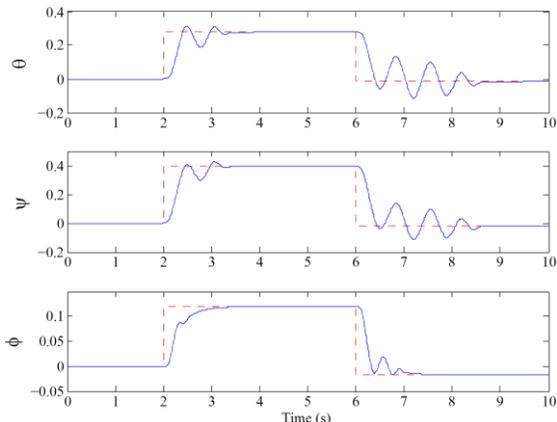

Fig. 5. Desired (blue line) and produced (red dashed line) virtual vectors

Fig. 3 shows the input vector and $w$ during the modified procedure.

## V. SIMULATION RESULTS

Consider the over-actuated Satellite Launch Vehicle (SLV) with a typical eight-actuator configuration described in [21]. Dynamic equations are:

$$\dot{\boldsymbol{x}} = A\boldsymbol{x} + B_u \boldsymbol{u} \qquad (68)$$
$$\boldsymbol{y} = C\boldsymbol{x}$$

where $\boldsymbol{x} = \begin{bmatrix} \theta & \dot{\theta} & \psi & \dot{\psi} & \phi & \dot{\phi} \end{bmatrix}^T$. The above equation could be written as [21, 22]

$$B_u = B_v B \qquad (69)$$

where $B_v \in \mathbb{R}^{n \times f}$ and $B \in \mathbb{R}^{f \times m}$. Hence, (68) could be simplified as

$$\dot{\boldsymbol{x}} = A\boldsymbol{x} + B_v \boldsymbol{v} \qquad (70)$$
$$\boldsymbol{v} = B\boldsymbol{u}$$

$\bar{\boldsymbol{u}}$ and $\underline{\boldsymbol{u}}$ represent the maximum and minimum physical limitations, respectively. Moreover, the actuators rate are considered as $\dot{\boldsymbol{u}}$ with the following limitations

$$\underline{\dot{\boldsymbol{u}}} \leq \dot{\boldsymbol{u}} \leq \bar{\dot{\boldsymbol{u}}} \qquad (71)$$

By considering $T$ as the sampling time, the overall position of constrains at time $t$ could be written as:

$$\boldsymbol{u}_{min} = \max(\underline{\boldsymbol{u}}, \boldsymbol{u}(t-T) + \underline{\dot{\boldsymbol{u}}}T)$$
$$\boldsymbol{u}_{max} = \min(\bar{\boldsymbol{u}}, \boldsymbol{u}(t-T) + \bar{\dot{\boldsymbol{u}}}T) \qquad (72)$$

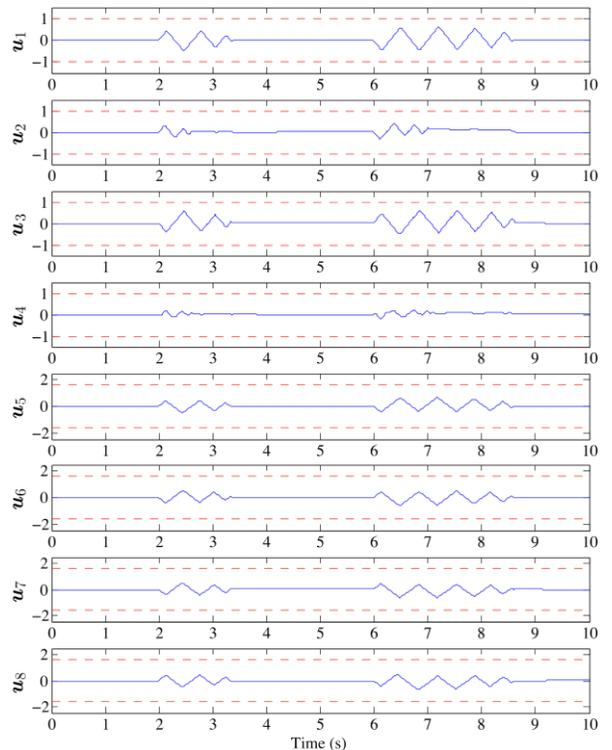

Fig. 6. Input vector position (blue line) within the maximum and minimum positions (red dashed line)

The numerical system data is provided in Appendix B. A linear quadratic set point control with state feedback is implemented as shown in fig. 4 [23]. $K$ is the LQR gain with the cost function defined as

$$J = \int_0^\infty (\boldsymbol{x}^T Q \boldsymbol{x} + \boldsymbol{u}^T R \boldsymbol{u}) dt \qquad (73)$$

where $Q$ and $R$ are the weighting matrices for the states and inputs. Note that there are linearly dependent columns in $B$ which make $\det(N_S N_S^T) = 0$. At times 2 and 6 seconds the following set points are implemented

$$\boldsymbol{r}_2 = [0.28 \quad 0.4 \quad 0.12]^T \frac{\pi}{180}$$

$$\boldsymbol{r}_6 = [-0.014 \quad -0.018 \quad -0.018]^T \frac{\pi}{180}$$

Fig. 5 shows the desired and produced values of the system states. It is shown that the controller and control allocation strategy could successfully follow the desired system states. Moreover, figures 6 and 7 illustrate the positions of the system inputs and their rates, respectively. It is obvious that the control allocation could hold the inputs within their limitations.

## VI. CONCLUSION

In this paper, a pseudo inverse based method is proposed to develop a new methodology which could solve the control allocation problem for the infeasible solution. Also, an appropriate method to find the matrix inverse of a singular matrix for the control allocation applications is proposed. The proposed method is based on the pseudo inverse, least squares, and singular value decomposition concepts. In contrast to the linear and quadratic programming methodologies, the number of iterations in this method is known. Numerical and simulation results are employed to show the effectiveness of the proposed method.

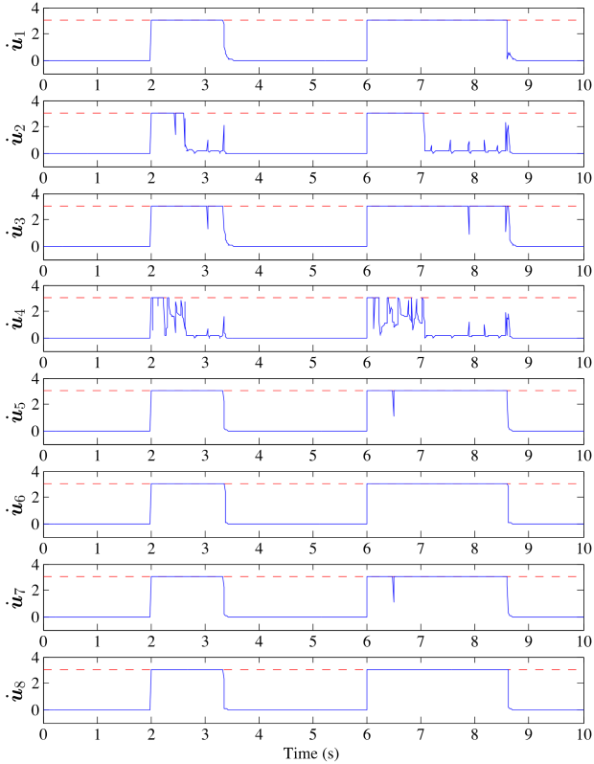

Fig. 7. Input vector rate position (blue line) within the maximum and minimum rate positions (red dashed line)

## APPENDIX A

**Lemma 11.** Assume that $N \in \mathbb{R}^{m \times n}$ where $m > n$. If some rows are removed, the resulted matrix is denoted by $N'$, and the following inequalities regarding its singular values hold:

$$\sigma' \leq \sigma$$
$$\underline{\sigma}' \leq \underline{\sigma} \quad (74)$$
$$\bar{\sigma}' \leq \bar{\sigma}$$

where $\sigma$, $\underline{\sigma}$, and $\bar{\sigma}$ are the gain, minimum and maximum singular values of $N$ and $\sigma'$, $\underline{\sigma}'$, and $\bar{\sigma}'$ are the gain, minimum and maximum singular values of $N'$, respectively.

*Proof:* Consider the following equation

$$\|r\|_2 = \|Nv_f\|_2 \quad (75)$$

so

$$\|r\|_2 = \sigma \|v_f\|_2 \quad (76)$$

Eliminate an arbitrary row from $N$. Hence

$$\|r'\|_2 = \|N' v_f\|_2 \quad (77)$$

and

$$\|r'\|_2 = \sigma' \|v_f\|_2 \quad (78)$$

It is obvious that $\|r'\|_2 \leq \|r\|_2$. By comparing (76) with (78), it yields

$$\sigma' \leq \sigma$$

The other inequalities can be similarly shown. ∎

## APPENDIX B

$$A = \begin{bmatrix} 1 & 1 & 0 & 0 & 0 & 0 \\ -0.2934 & -1 & 1.87e-5 & 0 & 0 & 0 \\ 0 & 0 & 1 & 1 & 0 & 0 \\ 2.71e-5 & 0 & -0.5621 & -1 & 0 & 0 \\ 0 & 0 & 0 & 0 & 1 & 1 \\ 5.71e-4 & 0 & 5.468e-4 & 0 & -1 & -1 \end{bmatrix}$$

$$B_v = \begin{bmatrix} 0 & 0 & 0 \\ 28.7539 & -0.0011 & -0.0075 \\ 0 & 0 & 0 \\ 0.0011 & -25.8651 & 0.0104 \\ 0 & 0 & 0 \\ -0.0232 & 0.0439 & 70.7638 \end{bmatrix}$$

$$B = \begin{bmatrix} 0.2985 & -0.2985 & -0.2985 & 0.2985 & 0.0779 & 0 & -0.0779 & 0 \\ -0.2985 & -0.2985 & 0.2985 & 0.2985 & 0 & 0.0779 & 0 & -0.0779 \\ -0.1618 & 0.1618 & -0.1618 & 0.1618 & 0 & 0.0211 & 0 & 0.0211 \end{bmatrix}$$

$$C = \begin{bmatrix} 1 & 0 & 0 & 0 & 0 & 0 \\ 0 & 0 & 1 & 0 & 0 & 0 \\ 0 & 0 & 0 & 0 & 1 & 0 \end{bmatrix}$$

$$\underline{u} = \begin{bmatrix} -1 & -1 & -1 & -1 & -1.6 & -1.6 & -1.6 & -1.6 \end{bmatrix}^T \frac{\pi}{180} rad$$

$$\bar{u} = \begin{bmatrix} 1 & 1 & 1 & 1 & 1.6 & 1.6 & 1.6 & 1.6 \end{bmatrix}^T \frac{\pi}{180} rad$$

$$\underline{\dot{u}} = \begin{bmatrix} -3 & -3 & -3 & -3 & -3 & -3 & -3 & -3 \end{bmatrix}^T \frac{\pi}{180} rad$$

$$\dot{\bar{u}} = \begin{bmatrix} 3 & 3 & 3 & 3 & 3 & 3 & 3 & 3 \end{bmatrix}^T \frac{\pi}{180} rad$$

$$F = \begin{bmatrix} -1 & 0 & 0 \\ 1 & 0 & 0 \\ 0 & -1 & 0 \\ 0 & 1 & 0 \\ 0 & 0 & -1 \\ 0 & 0 & 1 \end{bmatrix}$$

$$N = \begin{bmatrix} -2.4574e-2 & -7.4275e-9 & 0 \\ -3.8343e-9 & 1.693e-2 & 0 \\ -1.6126e-5 & -1.823e-5 & 0 \end{bmatrix}$$

$$K = \begin{bmatrix} 4.451 & 1.11 & 0 & 0 & -0.0007 & 0.0001 \\ 0 & 0 & -4.455 & -1.121 & 0.0024 & 0.0006 \\ 0.0007 & 0.0002 & 0.0025 & 0.0006 & 4.362 & 1.046 \end{bmatrix}$$

## REFERENCES
11